\documentclass[12pt,leqno]{article}
\usepackage{amsmath, amsthm, amsfonts, amssymb, color}
\usepackage{mathrsfs}
\usepackage{color}
\newcommand{\red}{\color{red}}

\theoremstyle{definition}

\newcommand{\scr}[1]{\mathscr #1}
\definecolor{wco}{rgb}{0.5,0.2,0.3}

\numberwithin{equation}{section} \theoremstyle{remark}

\newcommand{\ua}{\uparrow}

\title{{\bf Derivative Formula and Harnack Inequality  for Linear SDEs Driven by L\'evy   Processes}\footnote{Supported in
 part by  Lab. Math. Com. Sys., NNSFC(11131003), SRFDP and the Fundamental Research Funds for the Central Universities.}
}

\author{
{\bf Feng-Yu Wang}\\
\footnotesize{School of Mathematical Sciences,
Beijing Normal University, Beijing 100875, China}\\
\footnotesize{and}\\ \footnotesize{Department of Mathematics,
Swansea University, Singleton Park, SA2 8PP, UK}\\
\footnotesize{Email: wangfy@bnu.edu.cn; F.Y.Wang@swansea.ac.uk}}

\begin{document}

\def\R{\mathbb R}  \def\ff{\frac} \def\ss{\sqrt} \def\B{\mathbf
B}
\def\N{\mathbb N} \def\kk{\kappa} \def\m{{\bf m}}
\def\dd{\delta} \def\DD{\Delta} \def\vv{\varepsilon} \def\rr{\rho}
\def\<{\langle} \def\>{\rangle} \def\GG{\Gamma} \def\gg{\gamma}
  \def\nn{\nabla} \def\pp{\partial} \def\EE{\scr E}
\def\d{\text{\rm{d}}} \def\bb{\beta} \def\aa{\alpha} \def\D{\scr D}
  \def\si{\sigma} \def\ess{\text{\rm{ess}}}
\def\beg{\begin} \def\beq{\begin{equation}}  \def\F{\scr F}
\def\Ric{\text{\rm{Ric}}} \def\Hess{\text{\rm{Hess}}}
\def\e{\text{\rm{e}}} \def\ua{\underline a} \def\OO{\Omega}  \def\oo{\omega}
 \def\tt{\tilde} \def\Ric{\text{\rm{Ric}}}
\def\cut{\text{\rm{cut}}} \def\P{\mathbb P} \def\ifn{I_n(f^{\bigotimes n})}
\def\C{\scr C}      \def\aaa{\mathbf{r}}     \def\r{r}
\def\gap{\text{\rm{gap}}} \def\prr{\pi_{{\bf m},\varrho}}  \def\r{\mathbf r}
\def\Z{\mathbb Z} \def\vrr{\varrho} \def\ll{\lambda}
\def\L{\scr L}\def\Tt{\tt} \def\TT{\tt}\def\II{\mathbb I}
\def\i{{\rm in}}\def\Sect{{\rm Sect}}\def\E{\mathbb E} \def\H{\mathbb H}
\def\M{\scr M}\def\Q{\mathbb Q} \def\texto{\text{o}} \def\LL{\Lambda}
\def\Rank{{\rm Rank}} \def\B{\scr B} \def\i{{\rm i}} \def\HR{\hat{\R}^d}

\baselineskip 20pt
\maketitle
\begin{abstract} By using lower bound conditions of the L\'evy measure, derivative formulae  and  Harnack inequalities are derived for linear stochastic differential equations driven by L\'evy processes.  As applications, explicit gradient estimates and  heat kernel inequalities are presented. As  byproduct,  a new Girsanov theorem for L\'evy processes is derived.

\end{abstract} \noindent

 \noindent
Running title:  Derivative Formula and Harnack Inequality
 \vskip 2cm

\section{Introduction}

 The derivative formula   enables one to derive explicit gradient estimates; while the Harnack inequality has been applied to the study of heat kernel estimates, contractivity properties, transportation-const inequalities and properties of the invrainat probability measures, see e.g.
 \cite{W10, DRW09, WY10} and references within (see \S 4.2 for some general results).

Recall that Bismut's derivative formula of elliptic diffusion semigroups \cite{Bismut},  also known as      Bismut-Elworthy-Li formula  due to \cite{EL},  is a powerful tool for stochastic analysis on Riemannian manifolds  and has been   extended and applied to  SDEs (stochastic differential equations) driven by noises with a non-trivial Gaussian parts, see e.g. \cite{TA} and references within for the study of diffusion-jump processes. But, up to our knowledge,  explicit  derivative formula relying only on the L\'evy measure is not yet available.

On the other hand, by using  couplings constructed through Girsanov transforms, the dimension-free Harnack inequality, first introduced by the author in \cite{W97} for diffusion semigroups on manifolds, has been established and applied to various  SDEs and SPDEs driven by Gaussian noises, see \cite{ATW06, ATW09, DRW09, LW, Ouyang, ORW, RW, W07, W10, WWX10, WX10, WY10, Z}.
Since   arguments  used in these  references essentially relies on special properties of the Brownian motion, they do not apply to the jump setting. Therefore, it is in particular interesting to  built up a reasonable theory on derivative formula and Harnack inequality for pure jump processes.

In this paper, we aim to establish derivative formula and Harnack inequality for the semigroup associated to SDEs driven by L\'evy jump processes using lower bound conditions of the L\'evy measure. As observed in a recent paper \cite{W11b}, where the coupling property is confirmed for a class of linear SDEs driven by L\'evy processes,  the Mecke formula on the Poisson space will play an alternative role in the jump case to the Girsanov transform in the diffusion case. Indeed, with helps of this formula we will be able to establish explicit derivative formulae and Harnack inequalities for  a class of  jump processes (see Sections 3 and 4).

Before move on, let us introduce some recent results concerning regularity properties of the semigroup associated to the following linear SDE
\beq\label{E} \d X_t= A_tX_t\d t +\si_t \d L_t \end{equation} on $\R^d$, where $A,\si: [0,\infty)\to \R^d\otimes \R^d$ are measurable such that   $\si_s$ is invertible for  $s\ge 0 $ and $A,\si,\si^{-1}$ are locally bounded,      $L_t$  is the L\'evy process on $\R^d$ with L\'evy measure $\nu$  (see e.g. \cite{A, J}).  Let $P_t$ be the semigroup associated to (\ref{E}), i.e.
$$P_t f(x)=\E f(X_t^x),\ \ \ t\ge 0, x\in\R^d,  f\in \B_b(\R^d),$$ where $\B_b(\R^d)$ is the set of all bounded measurable functions on $\R^d$, and $X_t^x$ is the solution with initial data $x$. To formulate the solution, for any $s\ge 0$ let   $(T_{s,t})_{t\ge s}$ solve the equation on $\R^d\otimes \R^d$:
$$\ff{\d}{\d t}T_{s,t}=A_t T_{s,t},\ T_{s,s}=I.$$ Let $T_t= T_{0,t}, t\ge 0.$ We have $T_t=T_{s,t}T_s$ for $t\ge s\ge 0$ and
\beq\label{E1}X_t=T_tx +\int_0^tT_{s,t}\si_s \d L_s,\ \ t\ge 0.\end{equation}
By using lower bound conditions of the L\'evy measure $\nu$, the coupling property and gradient estimates have been derived in
\cite{W11a, W11b, BSW, SW, SSW}. Moreover, by using subordinations, the dimension-free Harnack inequality has been established in \cite{GRW}  for some jump processes in terms of known inequalities in the diffusion setting, where for the log-Harnack inequality the associated Bernstein function can be very broad but for the Harnack inequality with a power the function was assumed to have a growth stronger than $\ss r$.
When   $A_t$ and $\si_t$ are independent of $t$ and $\nu(\d z)\ge c |z|^{-(\d+\aa)}\d z$ for some constants $c>0$ and $\aa\in (0,2)$, i.e. the equation is time-homogenous with noise having an $\aa$-stable part,   a different version of  Harnack inequality was presented in    \cite[Theorem 1.1 and Corollary 1.3]{JW}:   for any $p\ge 1$ there exists a constant $C>0$ such that
\beq\label{JW} (P_t f(x+h))^p\le C P_t f^p(x) \Big(1+\ff{|h|}{(t\land 1)^{\ff 1 \aa}}\Big)^{p(d+\aa)},\ \ t>0, x,h\in\R^d\end{equation}  holds for positive $f\in \B_b(\R^d)$. Since \eqref{JW} allows  $p=1$  which is impossible even for the Brownian motion,  this inequality is somehow stronger than {\red general}  ones derived in the diffusion setting. On the other hand, however, this equality is not sharp for small distance since when $h=0$ it is worse than the classical Jensen inequality. In particular, \eqref{JW} does not imply the strong Feller property as the usual ones do (see Theorem 4.4(2) below).

The remainder of the paper is organized as follows. In the next section we present two lemmas, which will be used to establish derivative formulae and Harnack inequalities in Sections 3 and 4 respectively. In particular, a new Girsanov theorem is presented for L\'evy processes using the L\'evy measure is presented, which is interesting by itself.   With concrete lower bounds of the L\'evy measure, explicit gradient estimates and Harnack inequalities will also be addressed in Sections 3 and 4, which extend and improve the corresponding known results derived recently in \cite{W11a, JW, SSW}, see Corollaries \ref{C1.3} and \ref{C1.5} for details.

\section{Prelimilary}

Form now on, let $t>0$ be fixed, and let
$$ W_t=\big\{w: [0,t]\to \R^d \text{\ is\ right-continuous\ with \ left\ limits}\big\}$$ be the path space, which is a Polish space under the Skorokhod metric.  Let $L=(L_s)_{s\in [0,t]}$ be {\red a} L\'evy process with L\'evy measure $\nu$, possibly with a Gaussian part and a drift part. THen its distribution $\LL$ is a probability measure on $W_t$. For any $w\in W_t$ and $s\in [0,t]$, let
$\DD w_s= w_s-w_{s-}=w_s {\red -    \lim_{s'\uparrow s} w_{s'}}.$ Then
\beq\label{MM0}w(\cdot):= \sum_{s\in [0,t], \DD w_s\ne 0}\dd_{(z,s)}\end{equation}  is a $\si$-finite $\Z_+\cap\{\infty\}$-valued measure on $\R^d\times [0,t]$. For any non-negative function $g$ on $\R^d\times [0,t]$, let
$$w(g)=\int_{\R^d\times [0,t]}g(z,s)w(\d z,\d s)=\sum_{s\in [0,t], \DD w_s\ne 0} g(z,s).$$

Now, let $\nu\ge \nu_0$, where $\nu_0$ is another L\'evy measure. We may  write $L=L^1+L^0$, where $L^1$ and $L^0$ are two independent L\'evy processes
with L\'evy measure $\nu-\nu_0$ and $\nu_0$ respectively, and $\LL^0$ does not have Gaussian term.    L
et $\LL^1$ and $\LL^0$ be the distributions of $L^1$ and $L^0$ respectively. We have $\LL=\LL^1*\LL^0.$ In the sequel we will mainly use the $L^0$ part to
establish derivative formulae of $P_t$.

It is well known that   $\LL^0$ can be represented by
using the Poisson measure $\Pi_t$  with intensity
$$\mu_t(\d z,\d s) =1_{[0,t]}(s)\nu_0(\d z)\times \d s, $$ which is a probability measure on the
configuration space
\beg{equation*}\beg{split} \GG_t:=\Big\{\gg:=\sum_{i=1}^n \dd_{(z_i,s_i)}:\ & n\in \Z_+\cup\{\infty\}, z_i\in\HR, s_i\in [0,t],\\ &\gg(\{|z|\ge\vv\}\times [0,t])<\infty\ \text{for}\ \vv>0\Big\}\end{split}\end{equation*} equipped with the vague topology, where  $\HR =\R^d\setminus \{0\}.$
More precisely (see e.g. \cite[(4.2)]{Wu}),
\beq\label{2.1} \LL^0= \Pi_t\circ \phi^{-1}\end{equation} holds for
$$\phi(\gg) :=  Bt+ \int_{[0,t]\times\{|z|>1\}} z1_{[s,t]}\gg(\d s,\d z) +\int_{[0,t]\times \{0<|z|\le 1\}}z1_{[s,t]}(\gg-\mu_t)(\d s,\d z) ,\ \ \gg \in\GG_t,$$ where $B\in \R^d$ is a constant. Since $\mu_t$ is a L\'evy measure on $[0,t]\times \R^d$,   $\phi(\gg)\in W_t$  is well-defined   for $\Pi_t$-a.s. $\gg$. It is easy to see that
\beq\label{MMM} \phi(\gg-\dd{(z,s)})= \phi(\gg)- z1_{[s,t]},\ \ \text{for}\ \dd_{(z,s)}\le\gg,\end{equation} and by (\ref{MM0}),
\beq\label{MM2} \gg(\d z,\d s)= \phi(\gg)(\d z,\d s).\end{equation}
This and the Mecke formula for the Poisson measure imply the following lemma, which is crucial for our study.
\beg{lem} \label{L2.1}  For any $h\in L^1(W_t\times\R^d\times [0,t]; \LL^0\times \nu_0\times \d s)$,
\beq\label{F1}\beg{split} &\int_{W_t\times\R^d\times [0,t]}  h(w,z,s) \LL^0(\d w)\mu_t(\d z,\d s)\\
&=\int_{W_t}\LL^0(\d w)\int_{\R^d\times [0,t]} h(w-z1_{[s,t]},z, s)w(\d z,\d s).\end{split}\end{equation}
 Consequently, for $X_t^x$ solving $(\ref{E})$ with initial data $x$, \beq\label{F2}\beg{split} &\E\int_{\R^d\times [0,t]} f(X_t^x + T_{s,t}\si_s z)h(L^0, z,s)\mu_t(\d z,\d s)\\
&= \E\int_{\R^d\times [0,t]} f(X_t^x)h(L^0-z1_{[s,t]}, z,s)L^0(\d z,\d s). \end{split}\end{equation}
\end{lem}

\beg{proof} {\red By the Mecke formula \cite{M} (see \cite[Lemma 6.7]{R}), for any $F\in L^1(\Pi_t\times \mu_t)$ we have
$$\int_{\GG_t} \Pi_t(\d\gg)\int_{\R^d\times [0,t]}F(\gg,z,s)\mu_t(\d z,\d s) =\int_{\GG_t} \Pi_t(\d\gg) \int_{\R^d\times [0,t]} F(\gg-\dd_{(z,s)},z,s)\gg(\d z,\d s).$$}
Combining  this with (\ref{2.1}), (\ref{MMM}) and (\ref{MM2}), we obtain
\beg{equation*}\beg{split}
 &\int_{W_t\times\R^d\times [0,t]}  h(w,z,s) \LL^0(\d w)\mu_t(\d z,\d s)\\
&=\int_{\GG_t}\Pi_t(\d\gg)\int_{\R^d\times [0,t]}  h(\phi(\gg),z,s) \mu_t(\d z, \d s) \\
&= \int_{\GG_t} \Pi_t(\d \gg)\int_{\R^d\times [0,t]}  h(\phi(\gg-\dd_{(z,s)}),z,s) \gg(\d z, \d s)\\
&=\int_{\GG_t}\Pi_t(\d \gg)\int_{\R^d\times [0,t]}  h(\phi(\gg)-z1_{[s,t]}),z,s) \phi(\gg)(\d z, \d s)\\
&= \int_{W_t}\LL^0(\d w)\int_{\R^d\times [0,t]}  h(w-z1_{[s,t]},z,s) w(\d z, \d s).\end{split}\end{equation*}   Hence, (\ref{F1}) holds.

 Next, let $$\psi(w)= T_t x+ \int_0^t T_{s,t}\si_s\d w_s,$$ where the integral w.r.t. $\d w_s$ is the It\^o integral which is $\LL$-a.s. defined on $W_t$. By (\ref{E1}) we have
 $$f(X_t^x)= f\circ\psi(L^1+L^0),\ \ f(X_t^x + T_{s,t}\si_s z)= f\circ\psi(L^1+L^0+ z1_{[s,t]}).$$ Combining this with (\ref{F1})  and noting that $L^1$ and $L^0$ are independent with distributions $\LL^1$ and $\LL^0$ respectively, we obtain
 \beg{equation*}\beg{split} &\E\int_{\R^d\times [0,t]} f(X_t^x + T_{s,t}\si_sz)h(L^0, z,s)\mu_t(\d z,\d s)\\
 &= \int_{W_t}\LL^1(\d w^1) \int_{W_t\times\R^d\times [0,t]} f\circ\psi(w^1+w^0+z1_{[s,t]}) h(w^0, z,s) \LL^0(\d w^0)\mu_t(\d z,\d s) \\
 &= \int_{W_t\times W_t} \LL^1(\d w^1) \LL^0(\d w^0)\int_{\R^d\times [0,t]} f\circ\psi(w^1+w^0) h(w^0-z1_{[s,t]}, z,s) w^0(\d z,\d s)\\
 & = \E\int_{\R^d\times [0,t]} f(X_t^x)h(L^0-z1_{[s,t]}, z,s)L^0(\d z,\d s).\end{split}\end{equation*}
 \end{proof}

  As an application of Lemma \ref{L2.1}, we have the following Girsanov theorem, which might be interesting by itself.

 \beg{thm} \label{G} Let $G\ge 0$ be a measurable function on $W_t\times \R^d\times [0,t]$ such that $(\LL^0\times \mu_t)(G)=1.$  Let $(\xi,\tau)$ be a random variable
 on $\R^d\times [0,t]$ such that the distribution of $(L^0,\xi,\tau)$ is $G(w,z,s)\LL^0(\d w)\mu_t(\d z).$ Let
 $$g(z,s)=\int_{W_t} G(w,z,s)\LL^0(\d w),$$ which is the distribution density of $(\xi,\tau)$ w.r.t. $\mu_t$. If $\mu_t(g>0)=\infty$ and
 $$1_{\{G(w,z,s)>0\}}g(z,s)= g(z,s),\ \ (\LL^0\times \mu_t){\rm-a.e.}$$ Then the process $$L^0+ \xi1_{[\tau,t]}
 :=\big(L^0_s+\xi 1_{[\tau,t]}(s)\big)_{s\in [0,t]}$$ has distribution $\LL^0$ under the probability measure $\Q:=R\P$, where
 $$R= \ff {g(\xi,\tau)} {G(L^0,\xi,\tau)\{L^0(g)+g(\xi,\tau)\}}.$$\end{thm}

 \beg{proof} Since $G$ is the distribution density of $(L^0,\xi,\tau)$ w.r.t. $\LL^0\times\mu_t$, we have $G(L^0,\xi,\tau)>0$ a.s. Similarly,
 $g(\xi,\tau)>0$ a.s. as well. Moreover, it is easy to see that $\E L^0(g)=1$ so that $w(g)<\infty$. Therefore, $R\in (0,\infty)$ a.s.

 Now,  for any non-negative measurable function $F$ on $W_t$, applying (\ref{F1}) to
 $$h(w,z,s)= \ff{F(w,+z1_{[s,t]})}{g(z,s)+ w(g)}1_{\{G>0\}}(w,z,s),$$ which is finite since $\mu_t(g>0)=\infty$ implies that $w(g)>0$ holds  $\LL^0$-a.e.,
 and using $1_{\{G(w,z,s)>0\}}g(z,s)=g(z,s)$, we obtain
 \beg{equation*}\beg{split} &\E_{\Q} F(L^0+\xi1_{[\tau,t]}) = \E\ff{F(L^0+\xi1_{[\tau,t]})g(\xi,\tau)}{G(L^0,\xi,\tau)\{g(\xi,\tau)+L^0(g)\}}\\
 &= \int_{G(w,z,s)>0} \ff{F(w+ z1_{[s,t]}) G(w,z,s)g(z,s)}{G(w,z,s)\{g(z,s)+w(g)\big\}}\LL^0(\d w)\mu_t(\d z,\d s)\\
 &= \int_{W_t\times\R^d\times [0,t]} \ff{F(w+ z1_{[s,t]}) g(z,s)}{g(z,s)+w(g)}\LL^0(\d w)\mu_t(\d z,\d s)\\
 &=\int_{W_t}\LL^0(\d w) \int_{\R^d\times [0,t]} \ff{F(w)g(z,s)}{w(g)}w(\d z,\d s)= \int_{W_t} F(w)\LL^0(\d w).\end{split}\end{equation*} This completes the proof.\end{proof}

 A simple choice of $G$ in the above Theorem is that $G(w,z,s)=g(z,s),$ i.e. $(\xi,\tau)$ is independent of $L^0.$
To derive gradient estimate from Theorem \ref{T1.1} below, we need the following $\GG$-function:
$$\GG(r)= \int_0^\infty s^{r-1}{\red \e^{-s}}\d s,\ \ r>0.$$

 \beg{lem} \label{L2.2} Let $\LL^0$ be the distribution of a L\'evy process with L\'evy measure $\nu_0$ which is not necessarily absolutely continuous w.r.t. the Lebesgue measure. Let $\mu_t(\d z,\d s)=\nu_0(\d z)\times \d s$ on $\R^d\times [0,t].$ Then for any non-negative measurable function $g$ on $\R^d\times [0,t]$,
 $$\int_{W_t} \ff{\LL^0(\d w)}{w(g)^\theta}=\ff 1 {\GG(\theta)}\int_0^\infty r^{\theta -1} \exp\Big[-\mu_t(1-\e^{-rg})\Big]\d r,\ \ \theta>0.$$
 \end{lem}

 \beg{proof} Noting that
 $$\ff 1 {s^\theta}=\ff 1 {\GG(\theta)} \int_0^\infty r^{\theta-1}\e^{-s r}\d r,\ \  s>0,$$ it follows from (\ref{2.1}) that
 \beg{equation*}\beg{split}\int_{W_t}\ff{\LL^0(\d w)}{w(g)^{\theta}} &=\int_{\GG_t} \ff{\Pi_t(\d\gg)}{\gg(g)^{\theta}}
= \ff 1 {\GG(\theta)}\int_0^\infty r^{\theta-1} \d r \int_{\GG_t} \e^{-r\gg(g)}\Pi_t(\d \gg)\\
 &=\ff 1 {\GG(\theta)}\int_0^\infty r^{\theta -1} \exp\Big[-\mu_t(1-\e^{-rg})\Big]\d r.\end{split}\end{equation*}\end{proof}

\section{Derivative formula and gradient estimates}

To establish  a derivative formula for $P_t$, we need an absolutely continuous lower bound of $\nu$. Let
$$\nu(\d z)\ge \nu_0(\d z):=\rr_0(z)\d z$$ such that $\nu_0(\R^d)=\infty$. Recall that the infinity of $\nu$ is essential to ensure the strong Feller property of $P_t$, which is necessary for the differentiability of the semigroup (see \cite{PZ} and references within for criteria on the strong Feller property).  Thus, the assumption $\nu_0(\R^d)=\infty$ is reasonable in order to establish a derivative formula of $P_t$.

Let $L^0=(L^0_s)_{s\in [0,t]}$ be the   L\'evy process with L\'evy measure $\nu_0$, and let
$L^1=(L^1_s)_{s\in [0,t]}$ be the L\'evy processes with L\'evy measures  $\nu_1:=\nu-\nu_0$ independent of $L^0$, so that $L:= L^0+L^1$ is the L\'evy process with L\'evy measure $\nu$   introduced above. Let $\HR=\R^d\setminus \{0\}$ and let $\nu_0(g)$ be the integral of $g$ w.r.t. $\nu_0$.

\beg{thm}\label{T1.1} Let $\nu(\d z)\ge \rr_0(z)\d z$ for some
non-negative $\rr_0\in W^{1,1}_{loc}(\HR)$ such that $\nu_0(\d
z):=\rr_0(z)\d z$ is an infinite measure.   If there exists a non-negative
measurable function $g$ on  $\R^d\times [0,t]$ differentiable in
$z\in\R^d$  such that \beq\label{C1} \int_0^\infty
\exp\Big[-\mu_t(1-\e^{-rg})\Big]\d r+ \int_{\R^d\times [0,t]}
\big\{\rr_0g+\rr_0|\nn g|+g|\nn\rr_0|\big\}(z,s) \d z\d
s<\infty,\end{equation} where $\nn$ is the gradient in $z\in\R^d$,
then for any $f\in \B_b(\R^d)$,
\beq\label{D} \beg{split} &\nn P_t f(x)\\ &= -\E \int_{\R^d\times [0,t]} f(X_t^x +T_{s,t}\si_s z) \ff{(\si_s^{-1}T_s)^*\{L^0(g)\nn(\rr_0g)+g^2\nn\rr_0\}}{(L^0(g)+g)^2}(z,s)\d z\d s\\
&= \E \bigg[f(X_t^x)\int_{\R^d\times [0,t]}
\ff{(\si_s^{-1}T_s)^*\{\rr_0g\nn
g-L^0(g)\nn(\rr_0g)\}}{L^0(g)^2\rr_0}(z,s)L^0(\d z,\d s)\bigg],
\end{split}\end{equation}{\red where $$L^0(g):= \int_{\R^d\times [0,t]}g\,\d L^0= \sum_{s\in [0,t], \DD L_s^0\ne 0} g(\DD L_s^0,s).$$}\end{thm}

\beg{proof} Noting that the second equality in (\ref{D})   follows from the first and (\ref{F2}), we only need to prove the first formula.

(a) We first prove for the case where $g$ has a compact support $K$. Let $\LL =\LL^0*\LL^1$ be the distribution of $L$.  For $f\in \B_b(\R^d)$ and $\vv\in (0,1)$, let
$$h_\vv(w)= f\bigg(T_t(x+\vv h) +\int_0^t T_{s,t}\si_s\d w_s\bigg).$$   By (\ref{E1}) and noting that $L^0$ and $L^1$ are independent with distributions $\LL^0$ and $\LL^1$ respectively, we have
$f(X_t^{x+\vv h})= h_\vv(L^0+L^1)$ and
$$P_t f(x+\vv h) = \int_{W_t\times W_t} h_\vv(w^1+w^0)\LL^0(\d w^0)\LL^1(\d w^1).$$ Since $T_t=T_{s,t}T_s$ for $s\in [0,t]$, and since due to (\ref{C1}) and Lemma \ref{L2.2}
 $w(g)>0$ holds for $\LL^0$-a.s. $w\in W_t$,   this implies
\beg{equation*}\beg{split} & P_t f(x+\vv h) = \int_{W_t\times W_t} \LL^1(\d w^1)\LL^0(\d w^0) \int_{\R^d\times [0,t]}\ff{h_\vv(w^1+w^0)g(z,s)}{w^0(g)}w^0(\d z,\d s)\\
&= \int_{W_t\times W_t} \LL^1(\d w^1)\LL^0(\d w^0) \int_{\R^d\times [0,t]}\ff{h_0(w^1+w^0+\vv \si_s^{-1}T_s h1_{[s,t]})g(z,s)}{w^0(g)}w^0(\d z,\d s).
\end{split}\end{equation*} Combining this with (\ref{F1})  for $\LL^0$ in place of $\LL$ and
$$h(w,z,s):= \ff{h_0(w^1+w+(z+\vv \si_s^{-1}T_sh)1_{[s,t]})g(z,s)}{(w^0+z1_{[s,t]})(g)}, \ \ w\in W_t,$$ we arrive at
\beg{equation*}\beg{split} &P_t f(x+\vv h)\\
 &= \int_{W_t} \LL^1(\d w^1) \int_{W_t\times \R^d\times [0,t]}\ff{h_0(w^1+w^0+(z+\vv \si_s^{-1}T_sh)1_{[s,t]})g(z,s)}
{w^0(g)+g(z,s)}\LL^0(\d w^0)\mu_t(\d z,\d s)\\
&=\int_{W_t} \LL^1(\d w^1) \int_{W_t\times \R^d\times [0,t]}\ff{h_0(w^1+w^0+(z+\vv \si_s^{-1}T_sh)1_{[s,t]})(\rr_0g)(z,s)}
{w^0(g) +g(z,s)}\LL^0(\d w^0)\d z\d s.\end{split}\end{equation*} Using the integral transform $z\mapsto z-\vv \si_s^{-1}T_sh$, it follows that
\beq\label{**}\beg{split}&P_t f(x+\vv h)=\\
&\int_{W_t} \LL^1(\d w^1) \int_{W_t\times \R^d\times [0,t]}\ff{h_0(w^1+w^0+z1_{[s,t]})(\rr_0g)(z-\vv \si_s^{-1}T_sh,s)}
{w^0(g)+g(z-\vv\si_s^{-1}T_sh,s)}\LL^0(\d w^0)\d z\d s.\end{split}\end{equation} Therefore,
\beq\label{Dh}\beg{split}& \ff{P_tf(x+\vv h)-P_tf(x)}\vv\\
& =\int_{W_t} \LL^1(\d w^1) \int_{W_t\times \R^d\times [0,t]}h_0(w^1+w^0+z1_{[s,t]})\Phi_\vv(w^0,z,s)\LL^0(\d w^0)\d z\d s\end{split}\end{equation} holds for
$$\Phi_\vv(w^0,z,s):= \ff 1 \vv \bigg(\ff{g\rr_0}{w^0(g)+g}(z-\vv \si_s^{-1}T_sh,s)- \ff{g\rr_0}{w^0(g)+g}(z,s)\bigg).$$ Since
\beg{equation*}\beg{split} \lim_{\vv\to 0} \Phi_\vv(w,z,s)&=-\Big\<\nn\ff{\rr_0g}{w(g)+g}(z,s),  \si_s^{-1}T_s h\Big\>\\
&= -\Big\<(\si_s^{-1}T_s)^*
\ff{w(g)\nn(\rr_0g)+g^2\nn\rr_0}{(w(g)+g)^2}(z,s),   h\Big\>,\end{split}\end{equation*}
to derive the desired derivative formula by letting $\vv\to 0$, we need to make use of the dominated convergence theorem. Since $g$ has a compact support and $ \sup_{s\in [0,t]} |\si_s^{-1}T_sh|<\infty$, there is a compact set $K\subset \R^d$ such that supp\,$\Phi_\vv\subset W_t\times K\times [0,t]$ holds for all $\vv\in (0,1)$. Since   $\LL^0(\d w)\times \d z \times \d s$ is finite on $W_t\times K\times [0,t]$, it suffices to show that $\{\Phi_\vv\}_{\vv\in (0,1)} $ is uniformly integrable w.r.t. this measure. Noting that
$$\Big|\nn \ff{\rr_0g}{w(g)+g}\Big|\le \ff{|\nn(\rr_0g)|}{w(g)+g} + \ff{\rr_0g|\nn g|}{(w(g)+g)^2}\le \ff{2(\rr_0|\nn g|+g|\rr_0|)}{w(g)},$$ there exists a constant $c>0$ such that
we have
\beg{equation*}\beg{split}   |\Phi_\vv (w^0,z,s)|
  &=   \bigg|\ff 1 \vv\int_0^\vv \bigg\{\ff{\d}{\d r}\Big(\ff{\rr_0g}{w^0(g)+g}\Big)(z-r\si_s^{-1}T_sh,s)\bigg\}\d r\bigg| \\
 &\le \ff {c} {\vv w(g)} \int_0^\vv \big\{\rr_0|\nn g|+g|\nn\rr_0|\big\}(z-r\si_s^{-1}T_sh,s)\d r.\end{split}\end{equation*}
 By (\ref{C1}) and Lemma \ref{L2.2} we see that $\int_{W_t} \ff 1 {w(g)} \LL^0(\d w)<\infty$. So, it suffices to show that
 $$\Psi_\vv(z,s):=\ff 1 \vv\int_0^\vv \big\{\rr_0|\nn g|+g|\nn\rr_0|\big\}(z-r\si_s^{-1}T_sh,s)\d r,\ \ \vv\in (0,1)$$ is uniformly integrable
 w.r.t. $\d z\times \d s $ on $K\times [0,t].$
 Since the function $s\mapsto (s-R)^+$ is convex, by the Jensen inequality we have
 $$(\Psi_\vv-R)^+(z,s)\le \ff 1 \vv\int_0^\vv \big(\rr_0|\nn g|+g|\nn\rr_0|-R\big)^+(z-r\si_s^{-1}T_sh,s)\d r.$$ So,
\beg{equation*}\beg{split} &\int_{K\times [0,t]} (\Psi_\vv(z,s)-R)^+\d z\d s\\
&\le \ff 1\vv \int_{\R^d\times [0,t]\times [0,\vv]} \big(\rr_0|\nn g|+g|\nn\rr_0|-R\big)^+(z-r\si_s^{-1}T_sh,s)\d z\d s\d r\\
&= \int_{\R^d\times [0,t]} \big(\rr_0|\nn g|+g|\nn\rr_0|-R\big)^+(z,s)\d z\d s,\ \ \vv\in(0,1),\end{split}\end{equation*}
where the last step is due to the integral transform $z\mapsto z+ r\e^{-sA}h$ for the integral w.r.t. $\d z$. Combining this with (\ref{C1}) we see that
$$\lim_{R\to\infty}\sup_{\vv\in (0,1)}\int_{K\times [0,t]} (\Psi_\vv(z,s)-R)^+\d z\d s=0,$$ that is, $\{\Psi_\vv\}_{\vv\in (0,1)}$ is uniform integrable w.r.t. $\d z\times\d s$ on $K\times [0,t].$

(b) Let $g$ satisfy (\ref{C1}). For any $n\ge 1$, let $g_n(z,s)=g(z,s)\{1\land (1+n-|z|)^+\}$ which has a compact support. By (a) we have
\beq\label{B}  \nn P_t f(x) = -\E \int_{\R^d\times [0,t]} f(X_t^x +\e^{(t-s)A} z)\ff{(\si_s^{-1}T_s)^*\{L^0(g_n)\nn(g_n\rr_0) +g^2_n\nn\rr_0\}}{(L^0(g_n)+g_n )^2}(z,s)\d z\d s\end{equation} Let $c=\e^{t\|A\|}$. It is easy to see that
$$\bigg|\ff{(\si_s^{-1}T_s)^*\{L^0(g_n)\nn(g_n\rr_0) +g^2_n\nn\rr_0\}}{(L^0(g_n)+g_n )^2}\bigg|
\le \ff{c\{\rr_0|\nn g|+ g |\nn\rr_0|+\rr_0g\} }{L^0(g_1)},\ \ n\ge 1 $$ holds for some constant $c>0$. Then,  according to  (\ref{C1}), the desired formula follows from the   dominated convergence theorem by letting $n\to\infty$ in (\ref{B}), provided
\beq\label{B2}\int_{W_t} \ff {\LL^0(\d w)} {w(g_1)} <\infty.\end{equation}  By Lemma \ref{L2.2} and (\ref{C1}) we have
\beg{equation*}\beg{split} &\int_{W_t} \ff {\LL^0(\d w)} {w(g_1)}= \ff 1 {\GG(1)}\int_0^\infty  \exp\Big[-\mu_t(1-\e^{-rg_1})\Big]\d r\\
&\le \ff 1 {\GG(1)}\int_0^\infty   \exp\Big[t\nu_0(|z|>1)-\mu_t(1-\e^{-rg})\Big]\d r<\infty\end{split}\end{equation*} since $\nu_0(|z|\ge 1)\le \nu(|z|\ge 1)<\infty$ as $\nu$ is a L\'evy measure. Therefore, the proof is finished.
 \end{proof}

\beg{cor}\label{C1.2}   Let $\rr_0\in W_{loc}^{1,1}(\HR)$ be non-negative such that $\nu(\d z)\ge \rr_0(z)\d z$, and let $g$ be a non-negative measurable function on $\R^d\times [0,t]$ differentiable in the first variable such that $\mu_t(g)<\infty$.    Then for any $p\in (1,\infty]$ and $f\in \B_b(\R^d)$,
\beg{equation*}\beg{split} |\nn P_t f|\le  &(P_t |f|^p)^{\ff 1 p} \bigg(\ff 1 {\GG(1)}\int_0^\infty \exp\Big[-\mu_t(1-\e^{-rg})\Big]\d r\bigg)^{\ff {p-1}p}\\
&\times \bigg( \int_{\R^d\times [0,t]} \Big\{\|\si_s^{-1}T_s\|\big(g|\nn \log \rr_0|+|\nn \log(\rr_0g)|\big)(z,s)\Big\}^{\ff p{p-1}} g(z,s) \mu_t(\d z,\d s)\bigg)^{\ff {p-1}p}.\end{split}\end{equation*}\end{cor}

\beg{proof} Assume that the desired upper bound   is finite. Then (\ref{C1}) holds. On the other hand, according to (\ref{F2}), for any $ f\in \B_b(\R^d)$ we have
\beq\label{F3} P_t   f(x)= \E \int_{W_t\times\R^d\times [0,t]} f(X_t^x+T_{s,t} \si_sz)\ff{(g\rr_0)(z,s)}{L^0(g)+g(z,s)} \d z\d s.\end{equation} Combining this with the first formula in (\ref{D}) and using the H\"older inequality, we obtain
\beg{equation*}\beg{split} &\  |\nn P_t f| \\
&\le   \E\int_{\R^d\times [0,t]} f(X_t^x+T_{s,t}\si_sz)\|\si_s^{-1}T_s\|\Big\{\big(g|\nn \log\rr_0|
+|\nn \log (g\rr_0)| \big) \ff{(g\rr_0)}{L^0(g)+g}\Big\}(z,s)\d z\d s\\
&\le (P_t |f|^p)^{\ff 1 p} \bigg( \E\int_{\R^d\times [0,t]} \Big\{\|\si_s^{-1}T_s\|\big(g|\nn\log \rr_0|+|\nn \log(\rr_0g)|\big)\Big\}^{\ff p{p-1}}  \ff{g}{L^0(g)+g}\d \mu_t \bigg)^{\ff {p-1}p}. \end{split}\end{equation*} Then the desired gradient estimate follows from Lemma \ref{L2.2} by noting that  $\ff{g}{L^0(g)+g}\le \ff {g}{L^0(g)}$ and $\E \ff 1 {L^0(g)} =\int_{W_t} \ff 1 {w(g)}\LL^0(\d w).$\end{proof}

To illustrate Corollary \ref{C1.2}, we present explicit conditions on the lower bound of $\nu$ for the gradient estimate and Harnack inequality.
Comparing with results in \cite{W11a, SSW} where the uniform gradient estimates are derived,  in the following result $\rr_0$ is not necessary corresponding to a Bernstein function and  more general $L^p$ gradient estimates are also provided. For a $d\times d$ matrix $M$ and a constant
$\aa\in\R$, we write $M\le\aa I$ provided $\<Ma,a\>\le \aa|a|^2$ holds for all $a\in\R^d.$

\beg{cor}\label{C1.3} Let $A_s\le \aa I$ and $\|\si_s^{-1}\|\le \ll$ for some constants $\aa\in\R$ and $\ll>0.$ Let $\nu(\d z)\ge |z|^{-d} S(|z|^{-2})1_{\{s\le r_0\}}$ for some constant $r_0>0$ and positive function $S\in C^1([0,\infty))$ such that
\beq\label{C2} \limsup_{r\to\infty} \ff{|S'(r)|r}{S(r)}<\infty. \end{equation} For  $p>1$ and  $k>2+\ff{p}{p-1}$, let
$$\psi_k(r)= \ff{(1-\e^{-1})\kk(d)} {2^k}\int_{\ff{r_0} 2 \land r^{-1/k}}^{\ff{r_0}2} \ff{S(s^{-2})}s \d s,\ \ \ r>0,$$ where $\kk(d)$ is the area of the unit sphere in $\R^d$. If $\int_0^\infty \e^{-t\psi_k(r)}\d r<\infty$,
then there exists a constant $c >0$ such that
$$|\nn P_t f|\le c\Big(\ff{\e^{\aa t}-1}\aa\Big)^{\ff{p-1}p}   (P_t |f|^p)^{\ff 1 p}  \bigg( \int_0^\infty \e^{-t\psi_k(r)}\d r  \bigg)^{\ff{p-1}p} $$ holds for $f\in \B_b(\R^d)$,  In particular, if $S(r)= c_0\log^{\vv}(1+r) $ for some $c_0,\vv>0,$ then for any $p>1$ there exists a constant $c>0$ such that
$$|\nn P_t f|\le (P_t f^p)^{1/p} \exp\big[ c  (1\land t)^{-\ff 1 \vv}\big],\ \ t>0$$ holds for all positive $f\in \B_b(\R^d)$.
\end{cor}

\beg{proof}  Let $\rr_0(z)= |z|^{-d}S(|z|^{-2}) {(1- \ff{|z|}{r_0})^+}^k$ and $g(z,s)=g(z)=|z|^k.$ Obviously, $\int_0^\infty \e^{-t\psi_k(r)}\d r<\infty$
implies that $\int_{\R^d}\rr_0(z)\d z=\infty$ . Since $\|S'\|_\infty<\infty$ implies that $S(s^{-2})\le cs^{-2}$ for some constant $c>0$ and all $s\le r_0$, we have
$$(g\rr_0)(z)\le c|z|^{k-d-2}1_{\{|z|\le r_0\}},$$ so that $\mu_t(g)=t\int_{\R^d} (\rr_0g)(z)\d z<\infty.$
Next, it is easy to see from (\ref{C2}) that
$$\big\{g|\nn\log \rr_0|+|\nn\log (\rr_0g)|\big\}(z)\le \ff{c}{|z|(r_0-|z|)^+},\ \ |z|<r_0$$ holds for some constant $c>0$.  Thus,  there exists a constant $c>0$ such that
$$\Big\{\rr_0g\big(|g\nn\log \rr_0|^{\ff p{p-1}}+|\nn\log (\rr_0g)|^{\ff p {p-1}}\big)\Big\}(z)\le c|z|^{k-d-2-\ff{p}{p-1}} {(r_0-|z|)^+}^{k-\ff p {p-1}}$$ holds for some constant $c>0$. Since $\|\si_s^{-1}T_s\|\le \ll\e^{\aa s}$ and $k>2+\ff p{p-1},$ this implies
 \beq\label{AA}  \beg{split} &\int_{\R^d\times [0,t]} \Big\{\|\si_s^{-1}T_s\|\big(g|\nn \log \rr_0|+|\nn \log(\rr_0g)|\big)(z,s)\Big\}^{\ff p{p-1}} g(z,s) \mu_t(\d z,\d s)\\
 &\le  \ff{\ll(\e^{\aa t}-1)}\aa\int_{\R^d} \Big\{\rr_0g\big(|g\nn\log \rr_0|^{\ff p{p-1}}+|\nn\log (\rr_0g)|^{\ff p {p-1}}\big)\Big\}(z)\d z \le  \ff{c(\e^{\aa t}-1)}\aa \end{split}\end{equation}  for some constant $c>0$.
Next, for $r\ge (2/r_0)^k$ we have
\beg{equation*}\beg{split} \nu_0(1-\e^{-rg})&\ge \ff 1 {2^k} \int_{|z|\le r_0/2} |z|^{-d}S(|z|^{-2})(1-\e^{-r|z|^k})\d z \\
&=\ff{\kk(d)} {2^k} \int_0^{r_0/2} s^{-1} S(s^{-2}) (1-\e^{-rs^k})\d s\\
& \ge \ff{\kk(d)} {2^k} \int_{\ff{r_0} {2^k} \land r^{-1/k}}^{\ff{r_0}2} \ff{S(s^{-2})(1-\e^{-1})}{s}\d s=\psi_k(r).\end{split}\end{equation*}
  Combining this  with (\ref{AA}), we prove the first assertion by Corollary \ref{C1.2}.

Next, let $S(r)=c_0\log^{\vv}(1+r)$. By the semigroup property and the Jensen inequality, it suffices to prove the desired gradient estimate for $t\in (0,1].$ It is easy to see that
$$ t \psi_k(r) \ge c_1 t\log^{1+\vv}(1+r) -c_2t \ge 2\log(1+r) - c_3 t^{- 1/\vv}$$ holds for some constants $c_1,c_2,c_3>0$.
Then the desired gradient estimate follows from the first part of this Corollary.
  \end{proof}

Note that the second estimate in Corollary \ref{C1.3} improves and extends  \cite[Example 1.3]{W11a} to $L^p$ gradient estimate with better short time behavior. On the other hand, however, Corollary \ref{C1.3} does not provide sharp estimate for the $\aa$-stable case. In general, to drive sharper gradient estimates, it might be necessary to take $g$ depending also  on $s$.

\section{Harnack inequality and applications}

We first investigate the Harnack inequality with a power in the sense of \cite{W97} and the log-Harnack inequality introduced in \cite{W10, RW}, then present some applications of these inequalities in an abstract framework. Recently, these type of inequalities have been established in \cite{GRW}  for some  jump processes with using subordinations from diffusion processes and in \cite{JW} using heat kernel bounds of the $\aa$-stable process.

\subsection{Harnack inequality}

For positive measurable functions $\rr_0,g$ on $\HR$, let $\nu_0(\d z)=\rr_0(z)\d z$ and
$$\gg_{\rr_0,g}(\theta,t)=\ff 1 {\GG(\theta)} \int_0^\infty r^{\theta-1} \exp\big[-t\nu_0(1-\e^{-rg})\big]\d r,\ \ \theta,t>0. $$

\beg{thm}\label{T1.4} Let $\aa\in\R$ and $\ll\ge 0$ be such that $A_s\le \aa I$ and $\|\si_s^{-1}\|\le\ll$ for $s\in [0,1].$ Let $\rr_0\in W_{loc}^{1,1}(\HR)$ and $g\in W_{loc}^{1,1} (\R^d)$ be positive such that $\nu(\d z)\ge \rr_0(z)\d z$ and $\nu_0(g>0)=\infty.$ Then for any $p>1$ and positive $f\in \B_b(\R^d)$,
\beg{equation*}\beg{split} &\ff{(P_t f)^p(x+h)}{  P_t f^p(x)}\\
& \le \bigg\{  \int_{W_t\times \R^d\times [0,t]}
\Big(\ff{(\rr_0g)(z)} {w(g)+g(z)}\Big)^{\ff p {p-1}}
\Big(\ff{w(g)+g(z+\si_s^{-1}T_sh)}{(\rr_0g)(z+\si_s^{-1}T_sh)}\Big)^{\ff 1 {p-1}}\LL^0(\d w)\d z\d s\bigg\}^{p-1}\\
&\le \exp\Big[\|\nn\log (\rr_0g)\|_\infty \ll\e^{\aa  } |h|\Big]
 \bigg\{1+\big(\ll\|\nn g\|_\infty\e^{\aa  }|h|\big)^{\ff 1 {(p-1)\lor 1}} \gg_{\rr_0,g}
\Big(\ff 1 {p-1},t\land 1\Big)^{(p-1)\land 1}\bigg\}^{(p-1)\lor 1}\end{split}\end{equation*}  holds for $x,h\in\R^d$ and $t>0.$  \end{thm}

  \beg{proof} Since $\nu_0(g>0)=\infty$ implies $w(g)>0$ for $\LL^0$-a.e. $w$, the right-hand side of the first inequality makes sense (could be infinite). Let $g(z,s)=g(z).$ By (\ref{**}) and the H\"orlder inequality, we obtain
\beg{equation*}\beg{split} &(P_t f(x+h))^p \\
&=\bigg\{ \int_{W_t} \LL^1(\d w^1) \int_{W_t\times \R^d\times [0,t]}
\ff{h_0(w^1+w^0+z1_{[s,t]})(\rr_0g)(z- \si_s^{-1}T_sh)} {w^0(g)+g(z-\si_s^{-1}T_sh)}\LL^0(\d w^0)\d z\d s\bigg\}^p \\
&\le \bigg\{\int_{W_t} \LL^1(\d w^1) \int_{W_t\times \R^d\times [0,t]}
\ff{h_0^p(w^1+w^0+z1_{[s,t]})(\rr_0g)(z)} {w^0(g)+g(z)}\LL^0(\d w^0)\d z\d s\bigg\}\\
&\qquad\times \bigg\{  \int_{W_t\times \R^d\times [0,t]}
\Big(\ff{(\rr_0g)(z- \si_s^{-1}T_sh)} {w(g)+g(z-\si_s^{-1}T_sh)}\Big)^{\ff p {p-1}}
\Big(\ff{w(g)+g(z)}{(\rr_0g)(z)}\Big)^{\ff 1 {p-1}}\LL^0(\d w)\d z\d s\bigg\}^{p-1}\\
&=  P_t f^p(x) \bigg\{  \int_{W_t\times \R^d\times [0,t]}
\Big(\ff{(\rr_0g)(z)} {w(g)+g(z)}\Big)^{\ff p {p-1}}
\Big(\ff{w(g)+g(z+\si_s^{-1}T_sh)}{(\rr_0g)(z+\si_s^{-1}T_sh)}\Big)^{\ff 1 {p-1}}\LL^0(\d w)\d z\d s\bigg\}^{p-1},\end{split}\end{equation*}
where in the last step we have used the transform $z\mapsto z+\si_s^{-1}T_sh$ for the integral w.r.t. $\d z$. This proves the first inequality.

Next, due to the semigroup property and the Jensen inequality, for the second inequality it suffices to consider $t\in (0.1].$ Then
\beq\label{Int} (P_t f(x+h))^p\le P_t f^p(x) \bigg(\int_{W_t\times \R^d\times [0,t]}\ff{g(z)}{w(g)+g(z)}
  (B_1B_2)^{\ff 1 {p-1}} \LL^0(\d w)\nu_0(\d z)\d s\bigg)^{p-1}\end{equation} holds for
$$B_1=B_1(w,z,s):= \ff{w(g)+g(z+\si_s^{-1}T_sh)}{w(g)+g(z)},\ \ B_2=B_2(w,z,s):= \ff{(\rr_0g)(z)}{(\rr_0g)(z+\si_s^{-1}T_sh)}.$$
Since $t\in (0,1]$ and $\|\si_s^{-1}T_s\|\le \ll\e^\aa$ for $s\in (0,1]$, we have
\beg{equation*}\beg{split} B_1\le & 1+\ff{|g(z+\si_s^{-1}T_sh)-g(z)|}{w(g)+g(z)}\le 1+\ff{\ll\|\nn g\|_\infty \e^{\aa }|h|}{w(g)+g(z)},\\
B_2\le & \exp\Big[\|\nn\log (\rr_0g)\|_\infty \ll\e^{\aa }|h|\Big],\ \ \ s\in [0,t].\end{split}\end{equation*} Moreover,  due to   Lemma \ref{L2.1}
\beq\label{E1'}\beg{split} &\int_{W_t\times \R^d\times [0,t]}\Big(1+\ff c {w(g)+g(z)}\Big)^{\ff 1 {p-1}}\ff{g(z)}{w(g)+g(z)}\LL^0(\d w)\nu_0(\d z)\d s\\
&=\int_{W_t}  \Big(1+\ff c {w(g)}\Big)^{\ff 1 {p-1}}\LL^0(\d w)\int_{\R^d\times [0,t]}\ff{g(z)}{w(g)}\nu_0(\d z)\d s\\
& =\int_{W_t} \Big(1+\ff c {w(g)}\Big)^{\ff 1 {p-1}}\LL^0(\d w) \end{split}\end{equation}
holds for $c\ge 0$. So, it follows from (\ref{Int}) that
$$\ff{(P_tf(x+h))^p}{ P_t f^p(x) } \le  \exp\Big[\|\nn\log (\rr_0g)\|_\infty \ll\e^{\aa }
 |h|\Big]\bigg(\int_{W_t}\Big(1+\ff{\|\nn g\|_\infty\ll\e^{\aa }|h|}{w(g)}\Big)^{\ff 1 {p-1}}\LL^0(\d w)\bigg)^{p-1}.$$ This implies the second inequality  since
$$\bigg(\int_{W_t}\Big(1+\ff{c}{w(g)}\Big)^{\ff 1 {p-1}}\LL^0(\d w)\bigg)^{p-1}
 \le   \bigg\{1+c^{\ff 1 {(p-1)\lor 1}}\gg_{\rr_0,g}
\Big(\ff 1 {p-1},t\Big)^{(p-1)\land 1}\bigg\}^{(p-1)\lor 1}$$ holds for $c\ge 0$ according to Lemma \ref{L2.2} and the triangle inequality for the norm
$$\|F\|_{\ff 1 {p-1}}:= \bigg(\int_{W_t}|F|^{\ff 1 {p-1}}(w)\LL^0(\d w)\bigg)^{(p-1)\land 1},\ \ r>0.$$
\end{proof}

Next, we consider the log-Harnack inequality.

\beg{thm}\label{TLH} Let $\aa,\ll$ and $\rr_0,g$ be in Theorem $\ref{T1.4}$. For any positive $f\in \B_b(\R^d)$,
$$P_t \log f(x+h)\le \log P_t f(x) +\ll\e^\aa|h| \big(\|\nn\log (\rr_0g)\|_\infty +\|\nn g\|_\infty \gg_{\rr_0,g}(1,t\land 1)\big)$$ holds for $t>0$ and $x,h\in\R^d$. \end{thm}

\beg{proof} Again, due to the semigroup property and the Jensen inequality, it suffices to prove for $t\in (0,1].$  Let
$$\OO(\d w^1,\d w^0, \d z,\d s)=\ff{g(z)}{w(g)+g(z)}\LL^1(\d w^1)\LL^0(\d w^0)\nu_0(\d z)\d s,$$ which is a probability measure on $W_t\times W_t\times \R^d\times [0,t]$ according to (\ref{**}) for $g(z,s)=g(z), f=1 $(hence, $h_\vv=1$) and $\vv=0.$ Let
$$G(w,z,s) = \ff{(w(g)+g(z))(\rr_0g)(z-\si_s^{-1}T_sh)}{(w(g)+g(z-\si_s^{-1}T_sh))(\rr_0g)(z)},$$ which is a probability density w.r.t. $\OO$ by the same reason. Moreover, using $\log f$ to replace $f$ in (\ref{**}) with $\vv=1$, we have
$$P_t \log f(x+h) =\int_{W_t\times W_t\times\R^d\times [0,t]}(\log h_0)(w^1+w^0+z1_{[s,t]})G(w^0,z,s)\OO(\d w^1,\d w^0, \d z, \d s).$$ So, by the Young inequality (see  \cite[Lemma 2.4]{ATW09}) and (\ref{**}) with $\vv=0$, we obtain
\beg{equation*}\beg{split} &P_t\log f(x+h)  \le \log \int_{W_t\times W_t\times\R^d\times [0,t]} h_0(w^1+w^0+z1_{[s,t]})\OO(\d w^1,\d w^0, \d z,\d s)+ \OO(G\log G)\\
&=\log P_t f(x)+ \int_{W_t\times \R^d\times [0,t]} \Big\{\ff{(\rr_0 g)(z-\si_s^{-1}T_sh)}{w(g)+g(z-\si_s^{-1}T_sh)}\log G(w,z,s)\Big\}\LL^0(\d w)\d z\d s,\end{split}\end{equation*} where $\OO(G\log G)$ is the integral of $G\log G$ w.r.t. the probability measure $\OO$. Since for $t\in (0,1]$ one has
$$G(w,z,s)\le \exp\bigg[\bigg(\ff{\|\nn g\|_\infty}{w(g)}+\|\nn \log (\rr_0g)\|_\infty\bigg)\ll\e^{\aa }|h|\bigg],$$ and since (\ref{E1'}) and the integral transform $z\mapsto z+\si_s^{-1}T_sh$ imply that
\beg{equation*}\beg{split} &\int_{W_t\times\R^d\times [0,t]}\ff{(\rr_0g)(z-\si_s^{-1}T_sh)}{w(g)+g(z-\si_s^{-1}T_sh)} \LL^0(\d w)\d z\d s\\
&=\int_{W_t\times\R^d\times [0,t]}\ff{g(z)}{w(g)+g(z)} \LL^0(\d w)\nu_0(\d z)\d s=1,\end{split}\end{equation*} we conclude that
$$P_t \log f(x+h) \le \log P_t f(x) +\ll\e^\aa |h| \bigg(\|\nn\log (\rr_0g)\|_\infty +\|\nn g\|_\infty\int_{W_t}\ff 1 {w(g)}\LL^0(\d w)\bigg).$$ This  completes the  proof according to Lemma \ref{L2.2}.\end{proof}

Finally, we consider a specific situation for $\nu$ having an $\aa$-stable like lower bound. Comparing with the Harnack inequality (\ref{JW}) derived recently in \cite{JW},    our   result (\ref{Har}) is better for small time and small $|h|$, and we only need the specific lower bound in a neighborhood of $0$.

 \beg{cor}\label{C1.5} Let $A_s$ and $\|\si_s^{-1}\|$ be bounded above, and let $\nu(\d z)\ge h(|z|)\d z$  for some positive decreasing
 function $h\in C^1((0,\infty))$ such that

\beq\label{C2'}  \sup_{r>0} \ff{|h'(r)|}{h(r)+h(r)^2}<\infty. \end{equation}  Then for any $p>1$ there exist two constants $c_1,c_2>0$ such that for any positive $f\in \B_b(\R^d)$,
\beq\label{Har} (P_t f(x+h))^p \le   P_t f^p(x) \e^{c_2|h|} \Big(1+c_2|h|^{\ff 1 {(p-1)\lor 1}} \int_0^\infty r^{\ff{2-p}{p-1}}
 \e^{-c_1(t\land 1) r(h^{-1}(r))^d}\d r \Big)^{(p-1)\lor 1}\end{equation}  holds for $t>0, x,h\in \R^d.$ Moreover, there exist  constants $c_1,c_2>0$ such that
\beq\label{LHar} P_t\log f(x+h)\le \log P_t f(x) +c_2 |h|\int_0^\infty \e^{-c_1(t\land 1) r(h^{-1}(r))^d} \d r,\ \ x,h\in\R^d, t>0\end{equation}  holds for positive $f\in \B_b(\R^d).$ \end{cor}

 \beg{proof}  Obviously, it suffices to prove for $t\in (0,1].$  Let  $\rr_0(z)=h(|z|)$ and $g(z)= \ff 1 {1\lor\rr_0(z)}=\ff 1 {1\lor h(|z|)}.$
  Then it is easy to see from (\ref{C2'}) that $\|\nn\log (\rr_0 g)\|_\infty, \|\nn g\|_\infty<\infty.$ Moreover, since
  $$(\rr_0 g)(z)= h(|z|)\land 1 = 1,\ \text{if}\ g(z)\le 1,$$ for $r\ge 1$ we have
$$\nu_0(1-\e^{-r g})\ge \ff r 2 \nu_0(g1_{\{g\le r^{-1}\}})\ge \ff{ \kk(d)r}2 \int_0^{h^{-1}(r)}   s^{d-1}\d s\ge c_1 (h^{-1}(r))^d$$
for some constants $c_1>0$. Thus, for any $\theta>1$, there exists constants $c_2>0$ such that
$$\gg_{\rr_0, g}(\theta,t)\le c_2 \int_0^\infty r^{\theta-1} \exp\big[ -tc_1 r (h^{-1}(r))^d\big]\d r$$ holds
for $\theta=\ff p{p-1}$. Therefore, (\ref{Har}) and (\ref{LHar})  follow  from Theorems \ref{T1.4} and \ref{TLH} respectively.
\end{proof}
To illustrate the Corollary \ref{C1.5}, we consider $\nu(\d z)\ge bc_0 |z|^{-(d+\aa)}$ for some $c_0>0$ and $\aa\in (0,2).$ Letting $h(r)= c_0 r^{-(d+\aa)}$ we have
$$\int_0^\infty r^{\ff{2-p}{p-1}} \e^{-c_1 t r(h^{-1}(r))^d}\d r \le c' t^{-\ff{d+\aa}{\aa(p-1)}}$$ for some constant $c'>0$ and all $t\in (0,1].$
Therefore,
$$(P_t f(x+h))^p\le P_t f^p(x) \e^{c|h|} (t\land 1)^{-\ff{d+\aa}{\aa (p-1)}}$$ and
$$P_t\log f(x+h)\le \log P_t f(x) +\ff{c|h|}{(t\land 1)^{\ff{\aa+d}\aa}}$$ hold for all $t>0, x,h\in\R^d$ and positive $f\in \B_b(\R^d).$
One may also derive explicit Harnack and log-Harnack inequalities for the case that $\nu(\d z)\ge c_0 |z|^{-d} \log^\vv (1+ |z|^{-1})$ for some $c_0,\vv>0.$

\subsection{Applications} For applications of our results derived in this section, we introduce some applications of Harnack inequalities which are essentially organized or generalized from \cite{W10, WX10, WY10}.
As most results presented below are not yet well known, we include brief proofs for readers' convenience.

Let $E$ be a {\red topological} space with Borel $\si$-field $\B$, let $\B (E)$(resp. $\B_b(E), \B_b^+(E))$ denote  the set of all measurable (resp. bounded measurable, bounded non-negative measurable) functions on $E$, and let
 $C(E)$ (resp. $C_b(E), C_b^+(E))$ stands for the set of   continuous (resp. bounded continuous, bounded non-negative continuous) functions on $E$. We recall some notions which will be considered in this subsection.

\beg{defn} Let $\mu$ be a probability measure on $(E,\B)$, and let $P$ be a bounded linear operator on $\B_b(E)$. \beg{enumerate} \item[$(i)$]    $\mu$   is called   \emph{quasi-invariant}   of $P$, if
 $ \mu P$ is absolutely continuous w.r.t. $\mu$, where $(\mu P)(A):= \mu(P1_A),\ A\in\F$. If   $\mu P=\mu$ then $\mu$ is called an \emph{invariant probability measure} of $P$. \item[$(ii)$]  A   measurable function $p$ on $E^2$ is called the   \emph{kernel} of  $P$   w.r.t. $\mu$, if
  $$Pf= \int_E p(\cdot, y)f(y)\mu(\d y),\ \ \ f\in \B_b(E).$$ \item[$(iii)$]  $P$ is called a \emph{Feller } operator, if $PC_b(E)\subset C_b(E)$, while it is called a \emph{strong Feller} operator if $P\B_b(E)\subset C_b(E).$  \end{enumerate} \end{defn}

From now on, we let $P$ be a Markov operator given by $$Pf(x)= \int_E f(y)P(x,\d y),\ \ \ f\in \scr B_b(E), x\in E $$ for a transition probability measure $P(x,\d y)$.  We will consider the following general version of \emph{Harnack} type inequality for $P$:
\beq\label{0LP} \Phi(P f(x))\le  \{P\Phi(f)(y)\}\e^{ \Psi(x,y)},\ \ x,y\in E, f\in \B_b^+(E),\end{equation}
 where $\Phi$ is a non-negative function on $[0,\infty)$ and $\Psi$ is a  measurable non-negative function   on $E^2$. In particular,   the log-Harnack inequality and Harnack inequality with a power $p>1$ addressed above refer to   $\Phi(r)=\e^r$ and
 $\Phi(r)=r^p$ respectively.

\beg{thm}\label{0P} Let     $\mu$ be a quasi-invariant probability measure of $P$. Let $\Phi\in C^1([0,\infty))$ be an increasing function with $\Phi'(1)>0$ and
 $\Phi(\infty):=\lim_{r\to\infty}\Phi(r)=\infty$, such that $(\ref{0LP})$
 holds.
 \beg{enumerate}\item[$(1)$] For any $x,y\in E$, $P(x,\cdot)$ and $P(y,\cdot)$ are equivalent.
 \item[$(2)$] If $\lim_{y\to x} \{\Psi (x,y)+\Psi(y,x)\}=0$ holds for all $x\in E$, then $P$ is strong Feller. \item[$(3)$] $P$ has a kernel $p$ w.r.t. $\mu$, so that any invariant probability measure of $P$ is absolutely continuous w.r.t. $\mu$.
 \item[$(4)$]   $P$ has at most one invariant probability measure and if it has, the kernel of $P$ w.r.t. the invariant probability measure is strictly positive.
 \item[$(5)$] The kernel $p$ of $P$ w.r.t. $\mu$ satisfies
 $$\int_E p(x,\cdot)\Phi^{-1} \Big(\ff{p(x,\cdot)}{p(y,\cdot)}\Big)\d\mu\le \Phi^{-1}(\e^{\Psi(x,y)}),\ \ x,y\in E,$$
 where $\Phi^{-1}(\infty):=\infty$ by convention.  \item[$(6)$] If $r\Phi^{-1}(r)$ is convex for $r\ge 0$, then the kernel $p$ of $P$ w.r.t. $\mu$ satisfies
  $$\int_E p(x,\cdot)p(y,\cdot)\d\mu\ge \e^{-\Psi(x,y)},\ \ x,y\in E.$$\end{enumerate} \end{thm}

\beg{proof} We only need to prove (1) since other assertions are included in \cite[Proposition 3.1]{WY10}.  Let $A\in\B$ be such that $P(y,A)=0$. Applying (\ref{0LP}) to $n1_A$ we obtain
$$\Phi(nP1_A(x))\le \e^{\Psi(x,y)} P\Phi(n1_A)(y) = \e^{\Psi(x,y)} \Phi(0).$$ Since $\Phi(r)\to \infty$ as $r\to\infty$, letting $n\to\infty$ we conclude that $P(x,A)=P1_A(x)=0$. That is, $P(x,\cdot)$ is absolutely continuous w.r.t. $P(y,\cdot)$ and vice versa.
\end{proof}

By Theorem \ref{0P}(1), if (\ref{0LP}) holds then
$$p_{x,y}(z):= \ff{P(x,\d z)}{P(y,\d z)}$$ exists. We aim to describe this function using the Harnack inequality.
For simplicity, we only consider  the Harnack inequality with a power $p>1$
  \beq\label{0H3} (Pf(x))^p\le (Pf^p(y))\e^{\Psi(x,y)},\ \ \ x,y\in E,
  f\in \scr B_b^+(E)\end{equation} and the log-Harnack inequality
  \beq\label{0LH'}  P(\log f)(x)\le \log Pf(y)+ \Psi(x,y),\ \ \ x,y\in E, f\ge 1, f
\in \scr B_b(E).\end{equation}
The following result is organized from \cite[Section 2]{W10}.

  \beg{prp}\label{P2.4}   $(\ref{0H3})$ holds if and only if $p_{x,y}$ exists and satisfies
   \beq\label{0H4} P\big\{ p_{x,y}^{1/(p-1)}\big\}(x)\le
   \e^{\Psi(x,y)/(p-1)},\ \ \ x,y\in E;\end{equation}
while $(\ref{0LH'})$ holds if and only if   $p_{x,y}$ exists and satisfies
   \beq\label{0H4''} P \{\log p_{x,y}\}(x)\le
   \Psi(x,y),\ \ \ x,y\in E.\end{equation} \end{prp}

 Finally, we consider the hyperbounded property and the entropy-cost inequality implied by (\ref{0H3}) and (\ref{0LH'}). Let $P$ have an invariant probability measure $\mu$. Then
 $\|\cdot\|_{p\to q}$ stands for the operator norm from $L^p(\mu)$ to $L^q(\mu)$. Moreover,  for a non-negative measurable function $\Psi$ on $E\times E$, and for $\C(\nu,\mu)$   the class of all couplings of $\mu$ and $\nu$, let
 $$W_\Psi(\mu,\nu) =\inf_{\pi\in\C(\nu,\mu)} \int_{E\times E} \Psi(x,y)\pi(\d x,\d y)$$ be the transportation-cost from $\nu$ to $\mu$ induced by the cost-function $\Psi$. The following result can be deduced as in the proof of \cite[Corollary 1.2]{RW} and \cite[Section 2]{W01}.

 \beg{prp}\label{P2.5} Let $P$ have an invariant probability measure $\mu$.   \beg{enumerate}\item[$(1)$] (\ref{0H3}) implies
 $$\|P\|_{p\to \dd p} \le \int_E \ff {\mu(\d x)} {\big\{\int_E \exp[-\Psi(x,y)]\mu(\d y)\big\}^\dd},\ \ \dd>1.$$
 \item[$(2)$] Let $P^*$ be the adjoint operator of $P$ in $L^2(\mu)$. Then $(\ref{0LH'})$ implies
 $$\int_E (P^*f)\log P^*f \d \mu \le W_\Psi(f\mu,\mu),\ \ f\ge 0, \int_Ef\d \mu=1.$$
 \end{enumerate} \end{prp}

\paragraph{Acknowledgement.} The author would like to thank the referee and Dr. Jian Wang for corrections and useful comments.

\beg{thebibliography}{99}

\bibitem{A}  Applebaum, D. 2004. \emph{L\'evy Processes and Stochastic Calculus.} Cambridge University Press.

\bibitem{ATW06}   Arnaudon, M.,    Thalmaier, A., and   Wang, F.-Y. 2006.
 Harnack inequality and heat kernel estimates
  on manifolds with curvature unbounded below. \emph{Bull. Sci. Math.} 130: 223--233.

\bibitem{ATW09} Arnaudon, M.,    Thalmaier, A., and   Wang, F.-Y. 2009.
 Gradient estimates and Harnack inequalities on non-compact Riemannian manifolds.
    \emph{Stoch. Proc. Appl.} 119: 3653--3670.

  \bibitem{Bismut}  Bismut, J. M. 1984. \emph{Large Deviations and the
Malliavin Calculus.} Boston: Birkh\"auser, MA.

 \bibitem{BSW}   B\"ottcher, B.,   Schilling, R. L., and   Wang, J. 2011. Constructions of coupling processes for L\'evy processes.
 \emph{Stoch. Proc. Appl.}
121: 1201--1216.

\bibitem{DRW09}  Da Prato, G., R\"ockner, M., and  Wang, F.-Y. 2009. Singular stochastic equations on Hilbert
spaces: Harnack inequalities for their transition semigroups. \emph{J. Funct. Anal.} 257: 992--017.

\bibitem{EL}    Elworthy, K. D.,   and    Li, X.-M.  1994. Formulae for the
derivatives of heat semigroups. \emph{J. Funct. Anal.} 125:
252--286.

 \bibitem{GRW}  Gordina, M.,  R\"ockner, M., and  Wang, F.-Y. 2011. Harnack
inequalities for subordinate semigroups.  \emph{Potential Analysis} 34: 293--307.

\bibitem{GW11}  Guillin, A., and   Wang, F.-Y. 2012. Degenerate Fokker-Planck equations : Bismut formula, gradient estimate  and Harnack inequality.  \emph{J. Diff. Equations} 253: 20--40.

\bibitem{J}   Jacob, N. 2001. \emph{Pseudo Differential Operators and Markov
Processes  (Volume I),}  Imperial College Press, London.

\bibitem{LW}   Liu, W., and    Wang, F.-Y. 2008. Harnack inequality and strong Feller
  property for stochastic fast diffusion equations. \emph{J. Math. Anal. Appl.}
  342: 651--662.

\bibitem{M}  Mecke, J. 1967. Stationaire zuf\"allige Ma$\bb$e auf lokalkompakten abelschen Gruppen.
\emph{Z. Wahrsch. verw. Geb.} 9: 36--58.

\bibitem{Ouyang}   Ouyang, S.-X. 2011.  Harnack inequalities and applications for multivalued stochastic evolution equations.  \emph{Infin. Dimens. Anal. Quant. Probab.  Relat. Topics.} 14: 261--278.

\bibitem{ORW}    Ouyang, S.-X.,   R\"ockner, M., and  Wang, F.-Y. 2012. Harnack inequalities and applications for Ornstein-Uhlenbeck semigroups with jump.
 \emph{ Pot. Anal.}  36: 301--315.

\bibitem{PZ}  Priola, E., and    Zabczyk, J. 2009.  Densities for Ornstein-Uhlenbeck processes with jumps. \emph{ Bull. Lond. Math. Soc.} 41: 41--50.

\bibitem{R}  R\"ockner, M. 1998. Stochastic analysis on configuration spaces: basic ideas and recent results. In \emph{ New Directions in Dirichlet Forms,}
157--231. AMS/IP Stud. Math. 8, Amer. Math. Soc., Providence, RI.

\bibitem{RW03}  R\"ockner, M., and   Wang, F.-Y. 2003. Harnack and functional inequalities for generalized Mehler semigroups. \emph{J. Funct. Anal.} 203: 237--261.

\bibitem{RW}  R\"ockner, M., and Wang, F.-Y. 2010. Log-Harnack  inequality for stochastic differential equations in Hilbert spaces and its consequences. \emph{Infin. Dimens. Anal. Quant. Probab.  Relat. Topics.} 13: 27--37.

\bibitem{SW}    Schilling, R. L., and   Wang, J. 2011.  On the coupling property of L\'evy processes. \emph{Inst. Henri Poinc. Probab. Stat.} 47: 1147--1159.

\bibitem{SSW}   Schilling, R. L.,    Sztonyk, P., and   Wang, J. 2012. Coupling property and gradient estimates of L\'evy processes via the symbol.  \emph{Bernoulli.} 18: 1128--1149.

\bibitem{TA}  Takeuchi, A. 2010. Bismut-Elworthy-Li-Type formula for stochastic differential equations with jumps.
\emph{ J. Theory Probab.}  23: 576--604.

\bibitem{W97}   Wang, F.-Y. 1997. On estimation of the
logarithmic Sobolev constant and gradient estimates of heat
semigroups. \emph{ Probab. Theory Relat. Fields.} 108: 87--101.

\bibitem{W01} Wang, F.-Y. 2001.  Logarithmic Sobolev inequalities: conditions and counterexamples. \emph{J. Operator Theory.} 46: 189--197.

\bibitem{W07} Wang, F.-Y. 2007. Harnack inequality and applications
for stochastic generalized porous media equations. \emph{Ann. Probab.} 35:
1333--1350.

\bibitem{W10} Wang, F.-Y. 2010. Harnack inequalities on manifolds with boundary and applications. \emph{J. Math. Pures Appl.} 94: 304--321.

\bibitem{W11a} Wang, F.-Y. 2011.  Gradient estimate for Ornstein-Uhlenbeck jump processes.  \emph{Stoch. Proc. Appl.} 121: 466--478.

\bibitem{W11b} Wang, F.-Y. 2011. Coupling for   Ornstein-Uhlenbeck   Processes with Jumps.  \emph{Bernoulli.} 17: 1136--1158.

\bibitem{WWX10}    Wang, F.-Y.,   Wu, J.-L., and  Xu, L. 2011. Log-Harnack inequality for stochastic Burgers equations and applications.
\emph{J. Math. Anal. Appl.} 384: 151--159.

\bibitem{WX10}   Wang, F.-Y., and    Xu, L. 2013. Derivative formula and applications for hyperdissipative stochastic Navier-Stokes/Burgers equations. \emph{Infin. Dimens. Anal. Quant. Probab.  Relat. Topics.} to appear.

\bibitem{WY10}   Wang, F.-Y., and   Yuan, C. 2011. Harnack inequalities for functional SDEs with multiplicative noise and applications.   \emph{Stoch. Proc. Appl.}    121: 2692--2710.

\bibitem{JW}  Wang, J. 2011. Harnack inequalities for Ornstein-Uhlenbeck processes driven by L\'evy processes.  \emph{Statist. Probab. Lett.} 81:  1436--1444.

\bibitem{Wu}   Wu, L. 2000. A new modified logarithmic Sobolev inequality for Poisson point processes and several applications. \emph{Probab. Theory Relat. Fields.} 118: 427--438.

\bibitem{Z}  Zhang, T.-S. White noise driven SPDEs with reflection: strong Feller properties and Harnack inequalities. \emph{Potential Anal.}
 33: 137--151.
 \end{thebibliography}

\end{document}